\documentclass[conference]{IEEEtran}
\IEEEoverridecommandlockouts
\usepackage[usenames]{color}
\usepackage[normalem]{ulem}

\newcommand{\wubd}[1]{}

\newcommand{\wubn}[1]{}
\renewcommand{\wubn}[1]{{\color{red}****\textsc{waheed's note: #1}****}} 

\usepackage{cite}

\ifCLASSINFOpdf
\else
\fi
\usepackage{graphicx}
\usepackage{epsfig}
\usepackage{amssymb}
\usepackage{amsmath}
\usepackage{amsthm}
\usepackage{multirow}
\usepackage{algorithm}
\usepackage{algorithmic}
\usepackage{array}

\theoremstyle{plain}
\newtheorem{theorem}{Theorem}
\newtheorem{lemma}[theorem]{Lemma}

\theoremstyle{definition}
\newtheorem{definition}[theorem]{Definition}

\theoremstyle{remark}

\begin{document}
%
\title{Frame Coherence and Sparse Signal Processing}



%
\author{
\IEEEauthorblockN{
Dustin G. Mixon\IEEEauthorrefmark{1},
Waheed U. Bajwa\IEEEauthorrefmark{2},
Robert Calderbank\IEEEauthorrefmark{2}
}
\IEEEauthorblockA{\IEEEauthorrefmark{1}Program in Applied and Computational Mathematics, Princeton University, Princeton, New Jersey 08544}
\IEEEauthorblockA{\IEEEauthorrefmark{2}Department of Electrical and Computer Engineering, Duke University, Durham, North Carolina 27708}
\thanks{
This work was supported by the Office of Naval Research under Grant N00014-08-1-1110,
by the Air Force Office of Scientific Research under Grants FA9550-09-1-0551 and
FA 9550-09-1-0643 and by NSF under Grant DMS-0914892.
Mixon was supported by the A.B. Krongard Fellowship.
The views expressed in this article are those of the authors and do not reflect the official policy or position of the United States Air Force, Department of Defense, or the U.S. Government.
}
}


\maketitle

\begin{abstract}
The sparse signal processing literature often uses random sensing matrices to obtain performance guarantees.
Unfortunately, in the real world, sensing matrices do not always come from random processes.
It is therefore desirable to evaluate whether an arbitrary matrix, or frame, is suitable for sensing sparse signals.
To this end, the present paper investigates two parameters that measure the coherence of a frame: worst-case and average coherence.
We first provide several examples of frames that have small spectral norm, worst-case coherence, and average coherence.
Next, we present a new lower bound on worst-case coherence and compare it to the Welch bound.
Later, we propose an algorithm that decreases the average coherence of a frame without changing its spectral norm or worst-case coherence.
Finally, we use worst-case and average coherence, as opposed to the Restricted Isometry Property, to garner near-optimal probabilistic guarantees on both sparse signal detection and reconstruction in the presence of noise.
This contrasts with recent results that only guarantee noiseless signal recovery from arbitrary frames, and which further assume independence across the nonzero entries of the signal---in a sense, requiring small average coherence replaces the need for such an assumption.
\end{abstract}


%
\IEEEpeerreviewmaketitle

\section{Introduction}

Many classical applications, such as radar and error-correcting codes, make use of over-complete spanning systems \cite{strohmer:acha03}.
Oftentimes, we may view an over-complete spanning system as a \emph{frame}.
Take $F=\{f_i\}_{i\in\mathcal{I}}$ to be a collection of vectors in some separable Hilbert space $\mathcal{H}$.
Then $F$ is a frame if there exist \emph{frame bounds} $A$ and $B$ with $0<A\leq B<\infty$ such that
$A\|x\|^2\leq\sum_{i\in\mathcal{I}}|\langle x,f_i\rangle|^2\leq B\|x\|^2$ for every $x\in\mathcal{H}$.
When $A=B$, $F$ is called a \emph{tight frame}.
For finite-dimensional unit norm frames, where $\mathcal{I}=\{1,\ldots,N\}$, the \emph{worst-case coherence} is a useful parameter:
\begin{equation}
\label{eq.mu defn}
\mu_F:=\max_{\substack{i,j\in\{1,\ldots,N\}\\i\neq j}}|\langle f_i,f_j\rangle|.
\end{equation}
Note that orthonormal bases are tight frames with $A=B=1$ and have zero worst-case coherence.
In both ways, frames form a natural generalization of orthonormal bases.

In this paper, we only consider finite-dimensional frames.
Those not familiar with frame theory can simply view a finite-dimensional frame as an $M\times N$ matrix of rank $M$ whose columns are the frame elements.
With this view, the tightness condition is equivalent to having the spectral norm be as small as possible; for an $M\times N$ unit norm frame $F$, this equivalently means $\|F\|_2^2=\frac{N}{M}$.

Throughout the literature, applications require finite-dimensional frames that are nearly tight and have small worst-case coherence \cite{candes:annstat09,donoho:tit06b,HP03,mixon:icassp11,strohmer:acha03,tropp:tit04,tropp:acha08,zahedi:acc10}.
Among these, a foremost application is sparse signal processing, where frames of small spectral norm and/or small worst-case coherence are commonly used to analyze sparse signals \cite{candes:annstat09,donoho:tit06b,tropp:tit04,tropp:acha08,zahedi:acc10}.
Recently, \cite{bajwa:jcn10} introduced another notion of frame coherence called \emph{average coherence}:
\begin{equation}
\label{eq.nu defn}
\nu_F:=\tfrac{1}{N-1}\max_{i\in\{1,\ldots,N\}}\bigg|\sum_{\substack{j=1\\j\neq i}}^N\langle f_i,f_j\rangle\bigg|.
\end{equation}
Note that, in addition to having zero worst-case coherence, orthonormal bases also have zero average coherence.
It was established in \cite{bajwa:jcn10} that when $\nu_F$ is sufficiently smaller than $\mu_F$, a number of guarantees can be provided for sparse signal processing.
It is therefore evident from \cite{bajwa:jcn10,candes:annstat09,donoho:tit06b,HP03,mixon:icassp11,strohmer:acha03,tropp:tit04,tropp:acha08,zahedi:acc10} that there is a pressing need for nearly tight frames with small worst-case and average coherence, especially in the area of sparse signal processing.

This paper offers four main contributions in this regard.
First, we discuss three types of frames that exhibit small spectral norm, worst-case coherence, and average coherence: normalized Gaussian, random harmonic, and code-based frames.
With all three frame parameters provably small, these frames are guaranteed to perform well in relevant applications.
Second, performance in many applications is dictated by worst-case coherence \cite{candes:annstat09,donoho:tit06b,HP03,mixon:icassp11,strohmer:acha03,tropp:tit04,tropp:acha08,zahedi:acc10}.
It is therefore particularly important to understand which worst-case coherence values are achievable.
To this end, the Welch bound \cite{strohmer:acha03} is commonly used in the literature.
However, the Welch bound is only tight when the number of frame elements $N$ is less than the square of the spatial dimension $M$ \cite{strohmer:acha03}.
Another lower bound, given in \cite{MSEA03} and \cite{XZG05}, beats the Welch bound when there are more frame elements, but it is known to be loose for real frames \cite{CHS96}.
Given this context, our next contribution is a new lower bound on the worst-case coherence of real frames.
Our bound beats both the Welch bound and the bound in \cite{MSEA03} and \cite{XZG05} when the number of frame elements far exceeds the spatial dimension.
Third, since average coherence is new to the frame theory literature, we investigate how it relates to worst-case coherence and spectral norm.
In particular, we want average coherence to satisfy the following property, which is used in \cite{bajwa:jcn10} to provide various guarantees for sparse signal processing:
\begin{definition}
We say an $M\times N$ unit norm frame $F$ satisfies the \emph{Strong Coherence Property} if
\begin{equation*}
\mbox{(SCP-1)}~~~\mu_F\leq\tfrac{1}{164\log N}\qquad\mbox{and}\qquad\mbox{(SCP-2)}~~~\nu_F\leq\tfrac{\mu_F}{\sqrt{M}},
\end{equation*}
where $\mu_F$ and $\nu_F$ are given by \eqref{eq.mu defn} and \eqref{eq.nu defn}, respectively.
\end{definition}
Since average coherence is so new, there is currently no intuition as to when (SCP-2) is satisfied.
As a third contribution, this paper shows how to transform a frame that satisfies (SCP-1) into another frame with the same spectral norm and worst-case coherence that additionally satisfies (SCP-2).
Finally, this paper uses the Strong Coherence Property to provide new guarantees on both sparse signal detection and reconstruction in the presence of noise.
These guarantees are related to those in \cite{candes:annstat09,donoho:tit06b,tropp:acha08}, and we elaborate on this relationship in Section V.
In the interest of space, the proofs have been omitted throughout, but they can be found in \cite{bajwa:preprint11}.

\section{Frame constructions}

Many applications require nearly tight frames with small worst-case and average coherence.
In this section, we give three types of frames that satisfy these conditions.

\subsection{Normalized Gaussian frames}
Construct a matrix with independent, Gaussian distributed entries that have zero mean and unit variance.
By normalizing the columns, we get a matrix called a \emph{normalized Gaussian frame}.
This is perhaps the most widely studied type of frame in the signal processing and statistics literature.

To be clear, the term ``normalized'' is intended to distinguish the results presented here from results reported in earlier works, such as \cite{bajwa:jcn10,baraniuk:ca08,candes:tit05,wainwright:tit09}, which only ensure that the frame elements of Gaussian frames have unit norm in expectation.
In other words, normalized Gaussian frames are frames with individual frame elements independently and uniformly distributed on the unit hypersphere in $\mathbb{R}^M$.

That said, the following theorem characterizes the spectral norm and the worst-case and average coherence of normalized Gaussian frames.

\begin{theorem}[Geometry of normalized Gaussian frames]
\label{thm.normalized gaussian frames}
Build a real $M\times N$ frame $G$ by drawing entries independently at random from a Gaussian distribution of zero mean and unit variance.
Next, construct a normalized Gaussian frame $F$ by taking $\smash{f_n:=\frac{g_n}{\|g_n\|}}$ for every $n=1,\ldots,N$.
Provided $\smash{60\log{N}\leq M\leq\frac{N-1}{4\log{N}}}$, then the following inequalities simultaneously hold with probability exceeding $1 - 11N^{-1}$:
\begin{enumerate}
\item[(i)] $\mu_F \leq \frac{\sqrt{15\log{N}}}{\sqrt{M} - \sqrt{12\log{N}}}$,
\item[(ii)] $\nu_F \leq \frac{\sqrt{15\log{N}}}{M - \sqrt{12M\log{N}}}$,
\item[(iii)] $\|F\|_2 \leq \frac{\sqrt{M} + \sqrt{N} + \sqrt{2\log{N}}}{\sqrt{M - \sqrt{8M\log{N}}}}$.
\end{enumerate}
\end{theorem}

\subsection{Random harmonic frames}
Random harmonic frames, constructed by randomly selecting rows of a discrete Fourier transform (DFT) matrix and normalizing the resulting columns, have received considerable attention lately in the compressed sensing literature \cite{candes:tit06a,candes:tit06b,vershynin:cpam08}.
However, to the best of our knowledge, there is no result in the literature that shows that random harmonic frames have small worst-case coherence.
To fill this gap, the following theorem characterizes the spectral norm and the worst-case and average coherence of random harmonic frames.

\begin{theorem}[Geometry of random harmonic frames]
\label{thm.random harmonic frames}
Let $U$ be an $N\times N$ non-normalized discrete Fourier transform matrix, explicitly, $U_{k\ell}:= \mathrm{e}^{2\pi\mathrm{i}k\ell/N}$ for each $k,\ell=0,\ldots,N-1$.
Next, let $\{B_i\}_{i=1}^N$ be a collection of independent Bernoulli random variables with mean $\smash{\frac{M}{N}}$, and take $\mathcal{M}:=\{i:B_i=1\}$.
Finally, construct an $|\mathcal{M}|\times N$ harmonic frame $F$ by collecting rows of $U$ which correspond to indices in $\mathcal{M}$ and normalize the columns.
Then $F$ is a unit norm tight frame: $\smash{\|F\|_2^2=\frac{N}{|\mathcal{M}|}}$.
Furthermore, provided $\smash{16\log{N}\leq M\leq \frac{N}{3}}$, the following inequalities simultaneously hold with probability exceeding $1 - 4N^{-1} - N^{-2}$:
\begin{enumerate}
\item[(i)] $\frac{1}{2}M \leq |\mathcal{M}| \leq \frac{3}{2}M$,
\item[(ii)] $\nu_F\leq\frac{\mu_F}{\sqrt{|\mathcal{M}|}}$,
\item[(iii)] $\mu_F \leq \sqrt{\frac{118(N-M)\log{N}}{MN}}$.
\end{enumerate}
\end{theorem}

\subsection{Code-based frames}

Many structures in coding theory are also useful for constructing frames.
Here, we build frames from a code that originally emerged with Berlekamp in \cite{B70}, and found recent reincarnation with \cite{YG06}.
We build a $2^m\times2^{(t+1)m}$ frame, indexing rows by elements of $\mathbb{F}_{2^m}$ and indexing columns by $(t+1)$-tuples of elements from $\mathbb{F}_{2^m}$.
For $x\in\mathbb{F}_{2^m}$ and $\alpha\in\mathbb{F}_{2^m}^{t+1}$, the corresponding entry of the matrix $F$ is
\begin{equation}
\label{eq.matrix defn}
F_{x\alpha}=\tfrac{1}{\sqrt{2^{m}}}(-1)^{\mathrm{Tr}\big[\alpha_0x+\sum_{i=1}^t\alpha_ix^{2^i+1}\big]},
\end{equation}
where $\mathrm{Tr}:\mathbb{F}_{2^m}\rightarrow\mathbb{F}_2$ denotes the trace map, defined by $\mathrm{Tr}(z)=\sum_{i=0}^{m-1}z^{2^i}$.
The following theorem gives the spectral norm and the worst-case and average coherence of this frame.

\begin{theorem}[Geometry of code-based frames]
The $2^m\times2^{(t+1)m}$ frame defined by \eqref{eq.matrix defn} is unit norm and tight, i.e., $\|F\|_2^2=2^{tm}$, with worst-case coherence $\mu_F\leq \frac{1}{\sqrt{2^{m-2t-1}}}$ and average coherence $\smash{\nu_F\leq\frac{\mu_F}{\sqrt{2^{m}}}}$.
\end{theorem}

\section{Fundamental limits on worst-case coherence}

In many applications of frames, performance is dictated by worst-case coherence. 
It is therefore particularly important to understand which worst-case coherence values are achievable.
To this end, the following bound is commonly used in the literature:

\begin{theorem}[Welch bound \cite{strohmer:acha03}]
\label{thm.welch bound}
Every $M\times N$ unit norm frame $F$ has worst-case coherence $\mu_F\geq\sqrt{\tfrac{N-M}{M(N-1)}}$.
\end{theorem}

The Welch bound is not tight whenever $N>M^2$ \cite{strohmer:acha03}.
For this region, the following gives a better bound:

\begin{theorem}[\!\!\cite{MSEA03},\cite{XZG05}]
\label{thm.complex bound}
Every $M\times N$ unit norm frame $F$ has worst-case coherence $\mu_F\geq1-2N^{-1/(M-1)}$.
Taking $N=\Theta(a^M)$, this lower bound goes to $1-\frac{2}{a}$ as $M\rightarrow\infty$.
\end{theorem}

For many applications, it does not make sense to use a complex frame, but the bound in Theorem~\ref{thm.complex bound} is known to be loose for real frames \cite{CHS96}.
We therefore improve Theorem~\ref{thm.complex bound} for the case of real unit norm frames:

\begin{theorem}
\label{thm.bound}
Every real $M\times N$ unit norm frame $F$ has worst-case coherence
\begin{equation}
\label{eq.bound}
\mu_F\geq\cos\bigg[\pi\Big(\tfrac{M-1}{N\pi^{1/2}}~\tfrac{\Gamma(\frac{M-1}{2})}{\Gamma(\frac{M}{2})}\Big)^{\frac{1}{M-1}}\bigg].
\end{equation}
Furthermore, taking $N=\Theta(a^M)$, this lower bound goes to $\cos(\frac{\pi}{a})$ as $M\rightarrow\infty$.
\end{theorem}

\begin{figure}[t]
\centering
\includegraphics[width=\columnwidth]{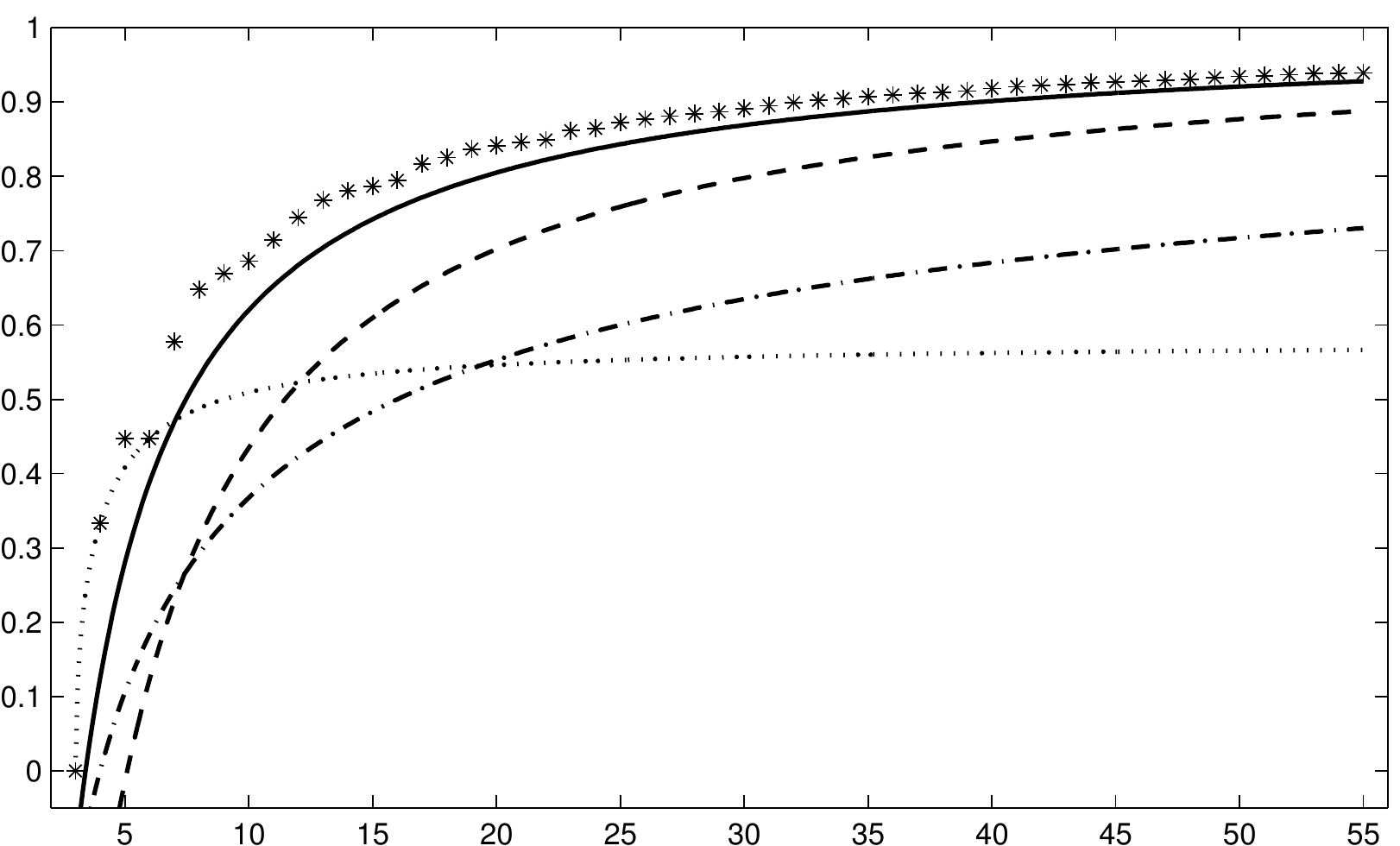}
\caption{Different bounds on worst-case coherence for $M=3$, $N=3,\ldots,55$.
Stars give numerically determined optimal worst-case coherence of $N$ real unit vectors, found in \cite{CHS96}.
Dotted curve gives Welch bound, dash-dotted curve gives bound from Theorem~\ref{thm.complex bound}, dashed curve gives general bound from Theorem~\ref{thm.bound}, and solid curve gives bound from Theorem~\ref{thm.3d points}.
 \label{figure}}
\end{figure}

In \cite{CHS96}, numerical results are given for $M=3$, and we compare these results to Theorems~\ref{thm.complex bound} and~\ref{thm.bound} in Figure~\ref{figure}.
Considering this figure, we note that the bound in Theorem~\ref{thm.complex bound} is inferior to the maximum of the Welch bound and the bound in Theorem~\ref{thm.bound}, at least when $M=3$.
This illustrates the degree to which Theorem~\ref{thm.bound} improves the bound in Theorem~\ref{thm.complex bound} for real frames.
In fact, since $\cos(\frac{\pi}{a})\geq 1-\frac{2}{a}$ for all $a\geq2$, the bound for real frames in Theorem~\ref{thm.bound} is asymptotically better than the bound for complex frames in Theorem~\ref{thm.complex bound}.
Moreover, for $M=2$, Theorem~\ref{thm.bound} says $\mu\geq\cos(\frac{\pi}{N})$, and \cite{BK06} proved this bound to be tight for every $N\geq2$.
For $M=3$, Theorem~\ref{thm.bound} can be further improved as follows:

\begin{theorem}
\label{thm.3d points}
Every real $3\times N$ unit norm frame $F$ has worst-case coherence $\mu_F\geq1-\frac{4}{N}+\frac{2}{N^2}$.
\end{theorem}

\section{Reducing average coherence}

In \cite{bajwa:jcn10}, average coherence is used to garner a number of guarantees on sparse signal processing.
Since average coherence is so new to the frame theory literature, this section will investigate how average coherence relates to worst-case coherence and the spectral norm.
We start with a definition:

\begin{definition}[Wiggling and flipping equivalent frames]
\label{def.flipping and wiggling}
We say the frames $F$ and $G$ are \emph{wiggling equivalent} if there exists a diagonal matrix $D$ of unimodular entries such that $G=FD$.
Furthermore, they are \emph{flipping equivalent} if $D$ is real, having only $\pm1$'s on the diagonal.
\end{definition}

The terms ``wiggling'' and ``flipping'' are inspired by the fact that individual frame elements of such equivalent frames are related by simple unitary operations.
Note that every frame with $N$ nonzero frame elements belongs to a flipping equivalence class of size $2^N$, while being wiggling equivalent to uncountably many frames.
The importance of this type of frame equivalence is, in part, due to the following lemma, which characterizes the shared geometry of wiggling equivalent frames:

\begin{lemma}[Geometry of wiggling equivalent frames]\label{lem:geom_eqframes}
Wiggling equivalence preserves the norms of frame elements, the worst-case coherence, and the spectral norm.
\end{lemma}

Now that we understand wiggling and flipping equivalence, we are ready for the main idea behind this section.
Suppose we are given a unit norm frame with acceptable spectral norm and worst-case coherence, but we also want the average coherence to satisfy (SCP-2).
Then by Lemma~\ref{lem:geom_eqframes}, all of the wiggling equivalent frames will also have acceptable spectral norm and worst-case coherence, and so it is reasonable to check these frames for good average coherence.
In fact, the following theorem guarantees that at least one of the flipping equivalent frames will have good average coherence, with only modest requirements on the original frame's redundancy.

\begin{theorem}[Frames with low average coherence]\label{thm:avc_rand}
Let $F$ be an $M\times N$ unit norm frame with $\smash{M < \frac{N-1}{4\log 4N}}$.
Then there exists a frame $G$ that is flipping equivalent to $F$ and satisfies $\smash{\nu_G\leq\frac{\mu_G}{\sqrt{M}}}$.
\end{theorem}

While Theorem~\ref{thm:avc_rand} guarantees the existence of a flipping equivalent frame with good average coherence, the result does not describe how to find it.
Certainly, one could check all $2^N$ frames in the flipping equivalence class, but such a procedure is computationally slow.
As an alternative, we propose a linear-time flipping algorithm (Algorithm~\ref{alg:flipping}).
The following theorem guarantees that linear-time flipping will produce a frame with good average coherence, but it requires the original frame's redundancy to be higher than what suffices in Theorem~\ref{thm:avc_rand}.

\begin{algorithm}[t]
\caption{Linear-time flipping}
\label{alg:flipping}
\textbf{Input:} An $M\times N$ unit norm frame $F$\\
\textbf{Output:} An $M\times N$ unit norm frame $G$ that is flipping equivalent to $F$
\begin{algorithmic}
\STATE $g_1\leftarrow f_1$ \hfill \COMMENT{Keep first frame element}
\FOR{$n=2$ to $N$}
\IF{$\|\sum_{i=1}^{n-1}g_i+f_n\|\leq\|\sum_{i=1}^{n-1}g_i-f_n\|$}
\STATE $g_n\leftarrow f_n$ \hfill \COMMENT{Keep frame element for shorter sum}
\ELSE
\STATE $g_n\leftarrow -f_n$ \hfill \COMMENT{Flip frame element for shorter sum}
\ENDIF
\ENDFOR
\end{algorithmic}
\end{algorithm}

\begin{theorem}
\label{thm.alg}
Suppose $N\geq M^2+3M+3$.
Then Algorithm~\ref{alg:flipping} outputs an $M\times N$ frame $G$ that is flipping equivalent to $F$ and satisfies $\nu_G\leq\frac{\mu_G}{\sqrt{M}}$.
\end{theorem}

As an example of how linear-time flipping improves average coherence, consider the following matrix:
\begin{equation*}
F:= \frac{1}{\sqrt{5}}\left[ \begin{array}{cccccccccc} +&+&+&+&-&+&+&+&+&-\\+&-&+&+&+&-&-&-&+&-\\+&+&+&+&+&+&+&+&-&+\\-&-&-&+&-&+&+&-&-&-\\-&+&+&-&-&+&-&-&-&- \end{array} \right].
\end{equation*}
Here, $\smash{\nu_F\approx0.3778>0.2683\approx\frac{\mu_F}{\sqrt{M}}}$.
Even though $N<M^2+3M+3$, we can run linear-time flipping to get the flipping pattern $D:=\mathrm{diag}(+-+--++-++)$.
Then $FD$ has average coherence $\smash{\nu_{FD}\approx0.1556<\frac{\mu_{F}}{\sqrt{M}}=\frac{\mu_{FD}}{\sqrt{M}}}$.
This example illustrates that the condition $N\geq M^2+3M+3$ in Theorem~\ref{thm.alg} is sufficient but not necessary.

\section{Near-optimal sparse signal processing without the Restricted Isometry Property}

Frames with small spectral norm, worst-case coherence, and/or average coherence have found use in recent years with applications involving sparse signals.
Donoho et al. used the worst-case coherence in \cite{donoho:tit06b} to provide uniform bounds on the signal and support recovery performance of combinatorial and convex optimization methods and greedy algorithms.
Later, Tropp \cite{tropp:acha08} and Cand\`{e}s and Plan \cite{candes:annstat09} used both the spectral norm and worst-case coherence to provide tighter bounds on the signal and support recovery performance of convex optimization methods for most support sets under the additional assumption that the sparse signals have independent nonzero entries with zero median.
Recently, Bajwa et al. \cite{bajwa:jcn10} made use of the spectral norm and both coherence parameters to report tighter bounds on the noisy model selection and noiseless signal recovery performance of an incredibly fast greedy algorithm called \emph{one-step thresholding (OST)} for most support sets and \emph{arbitrary} nonzero entries.
In this section, we discuss further implications of the spectral norm and worst-case and average coherence of frames in applications involving sparse signals.

\subsection{The Weak Restricted Isometry Property}
A common task in signal processing applications is to test whether a collection of measurements corresponds to mere noise \cite{kay:98b}.
For applications involving sparse signals, one can test measurements $y \in \mathbb{C}^M$ against the null hypothsis $H_0: y = e$ and alternative hypothesis $H_1: y = Fx+e$, where the entries of the noise vector $e\in \mathbb{C}^M$ are independent, identical zero-mean complex-Gaussian random variables and the signal $x\in\mathbb{C}^N$ is $K$-sparse.
The performance of such signal detection problems is directly proportional to the energy in $Fx$ \cite{davenport:jstsp10,haupt:icassp07,kay:98b}.
In particular, existing literature on the detection of sparse signals \cite{davenport:jstsp10,haupt:icassp07} leverages the fact that $\|Fx\|^2 \approx \|x\|^2$ when $F$ satisfies the Restricted Isometry Property (RIP) of order $K$.
In contrast, we now show that the Strong Coherence Property also guarantees $\|Fx\|^2 \approx \|x\|^2$ for most $K$-sparse vectors.
We start with a definition:

\begin{definition}
We say an $M\times N$ frame $F$ satisfies the \emph{$(K,\delta,p)$-Weak Restricted Isometry Property (Weak RIP)} if for every $K$-sparse vector $x \in \mathbb{C}^N$, a random permutation $y$ of $x$'s entries satisfies
\begin{equation*}
(1-\delta)\|y\|^2 \leq \|Fy\|^2 \leq (1+\delta)\|y\|^2
\end{equation*}
with probability exceeding $1-p$.
\end{definition}

We note the distinction between RIP and Weak RIP---Weak RIP requires that $F$ preserves the energy of \emph{most} sparse vectors.
Moreover, the manner in which we quantify ``most'' is important.
For each sparse vector, $F$ preserves the energy of most permutations of that vector, but for different sparse vectors, $F$ might not preserve the energy of permutations with the same support.
That is, unlike RIP, Weak RIP is \emph{not} a statement about the singular values of submatrices of $F$.
Certainly, matrices for which most submatrices are well-conditioned, such as those discussed in \cite{tropp:acha08}, will satisfy Weak RIP, but Weak RIP does not require this.
That said, the following theorem shows, in part, the significance of the Strong Coherence Property.

\begin{theorem}
Any $M\times N$ unit norm frame $F$ that satisfies the Strong Coherence Property
also satisfies the $(K,\delta,\frac{4K}{N^2})$-Weak Restricted Isometry Property
provided $N \geq 128$ and $2K\log{N} \leq \min\{\frac{\delta^2}{100\mu_F^2},M\}$.
\end{theorem}

\subsection{Reconstruction of sparse signals from noisy measurements}
Another common task in signal processing applications is to reconstruct a
$K$-sparse signal $x\in\mathbb{C}^N$ from a small collection of linear
measurements $y\in\mathbb{C}^M$. Recently, Tropp \cite{tropp:acha08} used both
the worst-case coherence and spectral norm of frames to find bounds on the
reconstruction performance of \emph{basis pursuit (BP)} \cite{donoho:siamjsc98}
for most support sets under the assumption that the nonzero entries of $x$ are
independent with zero median. In contrast, \cite{bajwa:jcn10} used the spectral
norm and worst-case and average coherence of frames to find bounds on the
reconstruction performance of OST for most support sets and \emph{arbitrary}
nonzero entries. However, both \cite{tropp:acha08} and \cite{bajwa:jcn10} limit
themselves to recovering $x$ in the absence of noise, corresponding to $y =
Fx$, a rather ideal scenario.

Our goal in this section is to provide guarantees for the reconstruction of
sparse signals from noisy measurements $y=Fx+e$, where the entries of the noise
vector $e\in \mathbb{C}^M$ are independent, identical complex-Gaussian random
variables with mean zero and variance $\sigma^2$. In particular, and in
contrast with \cite{donoho:tit06b}, our guarantees will hold for arbitrary
frames $F$ without requiring the signal's sparsity level to satisfy $K=O(\mu_F^{-1})$. 
The reconstruction algorithm that we analyze here is the OST
algorithm of \cite{bajwa:jcn10}, which is described in
Algorithm~\ref{alg:OST_recon}. The following theorem extends the analysis of
\cite{bajwa:jcn10} and shows that the OST algorithm leads to near-optimal
reconstruction error for large classes of sparse signals.

Before proceeding further, we first define some notation. We use
$\textsf{\textsc{snr}}:=\|x\|^2/\mathbb{E}[\|e\|^2]$ to denote the
\emph{signal-to-noise ratio} associated with the signal reconstruction problem.
Also, we use $\smash{\mathcal{T}_\sigma(t):=\{n: |x_n| >
\frac{2\sqrt{2}}{1-t}\sqrt{2 \sigma^2 \log{N}}\}}$ for any $t \in (0,1)$ to
denote the locations of all the entries of $x$ that, roughly speaking, lie
above the \emph{noise floor} $\sigma$. Finally, we use
$\smash{\mathcal{T}_\mu(t):=\{n: |x_n| > \frac{20}{t}\mu_F\|x\|\sqrt{2
\log{N}}\}}$ to denote the locations of entries of $x$ that, roughly speaking,
lie above the \emph{self-interference floor} $\mu_F\|x\|$.

\begin{algorithm}[t]
\caption{One-Step Thresholding (OST) \cite{bajwa:jcn10}}
\label{alg:OST_recon}
\textbf{Input:} An $M \times N$ unit norm frame $F$, a vector $y=Fx+e$, and a threshold $\lambda > 0$\\
\textbf{Output:} An estimate $\hat{\mathcal{K}}\subseteq\{1,\ldots,N\}$ of the support of $x$ and an estimate $\hat{x} \in \mathbb{C}^N$ of $x$
\begin{algorithmic}
\STATE $\hat{x} \leftarrow 0$ \hfill \COMMENT{Initialize}
\STATE $z \leftarrow F^* y$ \hfill \COMMENT{Form signal proxy}
\STATE $\hat{\mathcal{K}} \leftarrow \{n : |z_n| > \lambda\}$ \hfill \COMMENT{Select indices via OST}
\STATE $\hat{x}_{\hat{\mathcal{K}}} \leftarrow (F_{\hat{\mathcal{K}}})^\dagger y$ \hfill \COMMENT{Reconstruct signal via least-squares}
\end{algorithmic}
\end{algorithm}

\begin{theorem}[Reconstruction of sparse signals]
\label{thm:RSP}
Take an $M\times N$ unit norm frame $F$ which satisfies the
Strong Coherence Property, pick $t\in(0,1)$, and choose $\smash{\lambda =
\sqrt{2\sigma^2\log{N}}~\max \{\frac{10}{t}\mu_F\sqrt{M~\textsf{\textsc{snr}}},
\frac{\sqrt{2}}{1-t}\}}$. Further, suppose $x \in \mathbb{C}^N$ has support
$\mathcal{K}$ drawn uniformly at random from all possible $K$-subsets of
$\{1,\ldots,N\}$. Then provided
\begin{equation}
\label{thmeqn:RSP}
K \leq \tfrac{N}{c_1^2\|F\|_2^2\log{N}},
\end{equation}
Algorithm~\ref{alg:OST_recon} produces $\hat{\mathcal{K}}$ such that $\mathcal{T}_\sigma(t) \cap \mathcal{T}_\mu(t) \subseteq \hat{\mathcal{K}} \subseteq \mathcal{K}$ and $\hat{x}$ such that
\begin{equation}
\label{thmeqn:RSP_2}
\|x-\hat{x}\| \leq c_2 \sqrt{\sigma^2|\hat{\mathcal{K}}|\log{N}} + c_3\|x_{\mathcal{K} \setminus \hat{\mathcal{K}}}\|
\end{equation}
with probability exceeding $1 - 10N^{-1}$. Finally, defining $T:=|\mathcal{T}_\sigma(t) \cap \mathcal{T}_\mu(t)|$, we further have
\begin{equation}
\label{thmeqn:RSP_3}
\|x-\hat{x}\| \leq c_2 \sqrt{\sigma^2 K \log{N}} + c_3\|x - x_T\|
\end{equation}
in the same probability event.
Here, $c_1 = 37\mathrm{e}$, $c_2 = \frac{2}{1-\mathrm{e}^{-1/2}}$, and $\smash{c_3 = 1 + \frac{\mathrm{e}^{-1/2}}{1-\mathrm{e}^{-1/2}}}$ are numerical constants.
\end{theorem}

A few remarks are in order now for Theorem~\ref{thm:RSP}. First, if $F$
satisfies the Strong Coherence Property \emph{and} $F$ is nearly tight,
then OST handles sparsity that is almost linear in $M$: $K = O(M/\log{N})$ from
\eqref{thmeqn:RSP}. Second, the $\ell_2$ error associated with the OST
algorithm is the near-optimal (modulo the $\log$ factor) error of
$\sqrt{\sigma^2 K \log{N}}$ \emph{plus} the best $T$-term approximation error
caused by the inability of the OST algorithm to recover signal entries that are
smaller than $O(\mu_F\|x\|\sqrt{2 \log{N}})$. Nevertheless, it is easy to
convince oneself that such error is still near-optimal for large classes of
sparse signals. Consider, for example, the case where $\mu_F = O(1/\sqrt{M})$,
the magnitudes of $K/2$ nonzero entries of $x$ are some $\alpha =
\Omega(\sqrt{\sigma^2 \log{N}})$, while the magnitudes of the other $K/2$
nonzero entries are not necessarily same but scale as $O(\sqrt{\sigma^2
\log{N}})$. Then we have from Theorem~\ref{thm:RSP} that $\|x - x_T\| =
O(\sqrt{\sigma^2 K \log{N}})$, which leads to near-optimal $\ell_2$ error of
$\|x-\hat{x}\| = O(\sqrt{\sigma^2 K \log{N}})$. To the best of our knowledge,
this is the first result in the sparse signal processing literature that does not
require RIP and still provides near-optimal
reconstruction guarantees for such signals in the presence of noise, while using
either random or deterministic frames, even when $K = O(M/\log{N})$.

\end{document}